\documentclass{article}    

\usepackage{graphicx}          
\usepackage{epsfig}
\usepackage{graphicx}
\usepackage{makeidx}
\usepackage{multicol}
\usepackage{verbatim}
\usepackage{subfigure}
\usepackage{mathrsfs}
\usepackage{snapshot}

\usepackage{amsmath}

\usepackage{crop}
\crop

\usepackage{amssymb}
\usepackage{graphicx}

\def\E{\mathbb{E}}

\def\ba{\begin{array}}
\def\ea{\end{array}}
\def\bi{\begin{itemize}}
\def\ei{\end{itemize}}

\newtheorem{definition}{Definition}[section]

\newtheorem{theorem}{Theorem}[section]

\newtheorem{cond}{Condition}[section]
\begin{document}

\title{Recursive ECF identification of linear systems driven by L\'evy processes}
\date{}
\author{L\'aszl\'o Gerencs\'er and M\'at\'e M\'anfay}

\maketitle
In the literature the empirical characteristic function method is presented as an off-line identification method. While the results of the off-line methods are attractive, the proposed algorithms are ill-conditioned in many cases so that they requires special attention. As an alternative to the off-line method in this paper we propose and analyze on-line empirical characteristic function methods. Such recursive methods enables us to carry out real-time statistical analysis as new data points are processed instantly. In constructing these algorithms we follow the general framework proposed by Djereveckii and Fradkov , see \cite{djerev}, and Ljung, see \cite{ljung_scheme}. On-line methods are also used to complement a computationally expensive off-line identification method. Namely, it would be uneconomical to re-estimate $\theta^*$ using the off-line method when a new data point is received. Instead, we can argue that only a refinement of the estimate $\hat{\theta}_N$ should be computed using the newly received data point. This scenario not only shows a motivation behind the study of recursive algorithms but also suggests that it is reasonable to suppose that an initial guess of the parameter is close to $\theta^*.$

\section{General recursive estimation scheme}

We present a recursive estimation scheme within a general setup first formulated and solved for dynamical systems by Djereveckii and Fradkov in \cite{djerev} and Ljung in \cite{ljung_scheme}, hence the abbreviation DFL-scheme. Several recursive identification methods can be handled by this scheme, a nice summary of this can be found in \cite{soderstrom_rec}.
The basic building block of the scheme is the following parameter-dependent state-space equation:
\begin{equation}\label{eq:dynamic}
\overline{\xi}_{n+1}(x)=A(x)\overline{\xi}_n(x)+B(x)e_n,~\overline{\xi}_0(x)=0,
\end{equation}
where the parameter $x$ is an element of an open domain $D \subset \mathbb{R}^p.$ In the above so-called frozen parameter system $\overline{\xi} \in \mathbb{R}^r$ is a state-vector with possibly unobservable components and $e \in \mathbb{R}^m$ is an exogenous noise.
$x$ will be allowed to be time-varying taking values $(x_n)$ to be specified later. The next two conditions ensure the joint stability and the smoothness  of the matrices $A(x)$ and $B(x).$
\begin{cond}\label{cond:joint_stable}
The family of $r \times r$ matrices $\{ A(x), x \in D \subset \mathbb{R}^d\}$ is jointly stable, in the sense that there
exists a positive-definite $n \times n$ matrix $P,$ and a $\lambda$ with $0 < \lambda < 1$ such that
$$ A^T(x) P A(x) \leq \lambda P, $$
holds for all $x \in D.$
\end{cond}
\begin{cond}\label{cond:A_B_smooth}
$A(x)$ and $B(x)$ are continuously differentiable up to third order in $D.$
\end{cond}
To analyze recursive algorithms we require the driving noise process $(e_n)$ to be $L$-mixing, what is more we require that it is $L^+$-mixing in the sense defined below, defined in terms of the approximation error
$$
\gamma_q(\tau,e)=\sup_{n \geq \tau} \E^{1/q} \left[\left|e_n-\E\left[e_n|\mathscr{F}^+_{n-\tau}\right]\right|^q\right],
$$
see also the definition of $L$-mixing processes for example in \cite{gerencs_mixing}.
\begin{cond}\label{cond:e_n_strict}
$(e_n)$ is strictly stationary and it is also $L^+$-mixing with respect to a families of $\sigma$-algebras $(\mathscr{F}_n,\mathscr{F}_n^+)$ in the sense that for all integer $\tau \geq 1$ and $q \geq 1$ with some $c>0$ we have
$$\gamma_q(\tau,e)=O(\tau^{-1-c}).$$
\end{cond}
A variety of methods that analyze recursive methods is based on the idea of approximating $(x_n)$ using a trajectory of an ordinary differential equation (ODE).
In the process of developing the ODE method an often used assumption is that $e_n^2$ has some finite positive exponential moments. This leads to the definition of class $M^*.$
\begin{definition}
Let $(u_n), n \geq 0$ be a real-valued stochastic process. We say that $(u_n)$ is in class $M^*$ if for some $\varepsilon >0 $
$$
M^\varepsilon(u):=\sup_n \frac1\varepsilon \log \E \left[ e^{\varepsilon u_n} \right] < \infty.
$$
\end{definition}
\begin{cond}\label{cond:e_n}
$(e_n^2)$ is in $M^*.$
\end{cond}
Let $Q: \mathbb{R}^r \times \mathbb{R}^d \rightarrow \mathbb{R}^d$ denote a function such that it is bounded by some polynomial of $\overline{\xi}$ and the same holds for the derivatives of $Q$ up to order three. In many standard identification method $Q$ is a quadratic-form in $\overline{\xi}$, but for the present application this will not hold. Define
$$
F(x)=\lim_{n \rightarrow \infty} \E\left[Q(\overline{\xi}_n(x),x)\right].
$$
Now we are ready to formulate the abstract estimation problem related to the DFL-scheme: solve for $x$ the non-linear algebraic equation
\begin{equation}\label{eq:G=0}
F(x)=0.
\end{equation}
Without loss of generality we may assume that $x=x^*=0$ is a solution. We may assume that $D_0 \subset D$ is a compact domain such that $x^* \in D_0.$
Suppose that we are given an initial estimate of $x^*,$ say $x_0.$ Then the tentative recursion corresponding to the DFL scheme is given by
\begin{align}
\xi_{n+1} &= A(x_n)\xi_n + B(x_n)e_n,  \quad \xi_0=0. \label{ljung_rec_elso}\\
x_{n+1}&=x_n+\frac{1}{n+1}Q(\xi_{n+1},x_n) \quad x_0 \in D_0 \label{ljung_rec},
\end{align}
where $(x_n)$ denotes the sequence of generated estimates.
Typically the initial estimate $x_0$ is close to $x^*.$ A controversial issue is the problem of keeping $(x_n)$ in the domain $D_0.$ In order to guarantee this a resetting mechanism is introduced.
To make this modified recursion formal we denote the value of $x$ computed at time $n+1$ using (\ref{ljung_rec}) by $x_{n+1-}$ and define
$$
x_{n+1} =
  \begin{cases}
   x_{n+1-} & \text{if } x_{n+1-} \in D_0 \\
   x_0       & \text{if } x_{n+1-} \in D_0^c.
  \end{cases}
$$
That is, if $x_{n+1-}$ leaves the domain $D_0$ then a resetting is applied. This event is denoted by $B_{n+1}=\{ \omega | x_{n+1-} \in D_0^c\}.$ Hence
$(\ref{ljung_rec})$ is replaced by
\begin{equation}
x_{n+1}=x_{n}+(1-\boldsymbol{1}_{B_{n+1}})\frac{1}{n+1}Q(\xi_{n+1},x_n)+\boldsymbol{1}_{B_{n+1}}(x_0-x_n),
\end{equation}
where $\boldsymbol{1}_B$ is the indicator function of the event $B \subset \Omega$.

Now we define the differential equation, the solution trajectories of which reflect the pattern of behaviour of the sequence $(x_n).$ This so-called associated ordinary differential equation (ODE) is defined by
\begin{equation}\label{associated_od}
\dot{y}_t=\frac1t F(y_t),\quad y_s=\zeta,
\end{equation}
for $t\leq s \leq 1.$
Alternatively, the associated ODE can be also defined as
$$
\dot{y}_t=F(y_t),
$$
we allow this ambiguity in the definition of the ODE.
Note that since $F$ is well defined in $D$ and has continuous derivatives up to third order, the latter differential equation has a unique solution $y(t,s,\zeta)$ in some interval for $t$. 
A variety of convergence results is based on the stability of the above ODE (\ref{associated_od}), see \cite{soderstrom_rec}, \cite{chen_guo}, \cite{kusher} and \cite{benveniste2012adaptive}. For our application the stability of the associated ODE is specified by the next condition, which can be found in \cite{gerencser_recursive}.
\begin{cond}\label{cond:od_rate}
Let $D_0 \subset D$ be the compact truncation domain such that $x^* \in {\rm int} D_0.$ Assume that there exists a compact convex set $D'_{0}$ such that $D_0 \subset D'_0 \subset D$ and for all $t \geq s \geq 1$ we have
$$
y(t,s,\zeta) \in D'_{0} \text{~~for~~} \zeta \in D_0 \text{~~and~~} y(t,s,\zeta) \in D \text{~~for~~} \zeta \in D'_{0}.
$$
In addition $\lim_{t \rightarrow \infty}y(t,s,\zeta)=x^*$ for $\zeta \in D$ and
$$
\left\| \frac{\partial}{\partial \zeta} y(t,s,\zeta)\right\| \leq C(s/t)^{\alpha}
$$
with some $C>1,\alpha>0$ for all $\zeta \in D'_{0}$ and $t \geq s \geq 1.$ We have an initial estimate $\zeta=y_1=x_0$ such that for all $t \geq s \geq 1$ we have $y(t,s,\zeta) \in {\rm int} D_0.$ Finally, for the star-like closure
$$
D_0^*=\left\{y \left| y=x^*+\lambda(x-x^*), 0 \leq \lambda \leq 1, x \in D_0 \right. \right\}
$$ of the set $D_0$ we have $D_0^* \subset D.$
\end{cond}
The asymptotic covariance matrix of the estimates will be closely related to that of the averaged correction terms. Hence define
$$H(n,x,\omega):=Q(\overline{\xi}^{{\rm (s)}}_n(x),x)$$ and
the matrix $P^*$ in terms of $H$ as
\begin{equation}\label{eq:P^*}
P^*=\sum_{m=-\infty}^{\infty} \E \left[ H(m,x^*,\omega)H^*(0,x^*,\omega)\right].
\end{equation}
The following result from \cite{gerencser_recursive} states that the above recursion indeed defines a sequence of $x_n$-s that converges to the solution of the equation $F(x)=0$ and the rate of convergence of the moments of the error is also given. What is more the result also gives it asymptotic covariance matrix of the estimate.
\begin{theorem}\label{thm:rateconv_rec}
Assume that Conditions \ref{cond:joint_stable}-\ref{cond:e_n} are satisfied and further assume that
the differential equation (\ref{associated_od}) satisfies Condition \ref{cond:od_rate} with $\alpha > 1/2,$ then we have $x_N=O_M(N^{-1/2}).$ Moreover, the asymptotic covariance matrix of the error process $x_N-x^*,$ defined by
$$
\Sigma_{xx}=\lim_{N \rightarrow \infty} N \E\left[(x_N-x^*)(x_N-x^*)^* \right],
$$
exists and it satisfies the Lyapunov-equation
$$(A^*+I/2)\Sigma^*_{xx}+\Sigma_{xx}(A^*+I/2)^*+P^*=0,$$
where $A^*=F_{x}(x^*).$
\end{theorem}
An exciting special case is when the variable $x$ can be split as $x=(x_1,x_2)$ so that the recursive estimation method is partially stochastic Newton w.r.t $x_1$, meaning that
the Jacobian matrix of the r.h.s. of the corresponding associated ODE at $x=x^*$ is of the form
$$
\left(
  \begin{array}{cc}
    -I & 0 \\
    J_{2,1} & J_{2,2} \\
  \end{array}
\right).
$$
Then using simple linear algebra and Theorem \ref{thm:rateconv_rec} we conclude the following corollary.
{\bf Corrolary 1}
Assume that the conditions of Theorem \ref{thm:rateconv_rec} hold. Assume further that we can split $x$ as $x=(x^1,x^2)$ so that the recursive estimation method is a partially stochastic Newton method with respect to $x^1.$ Then the asymptotic covariance matrix of the recursive estimate $x^1_N$ equals to $P^*_{1,1},$ which is the corresponding block of $P^*$ defined in (\ref{eq:P^*}).


\section{Recursive ECF for i.i.d. sample}\label{sec:rec_noise}

The DFL-scheme provides a solution for the problem of estimating the parameters of a distribution or a regression function using i.i.d. samples by simply choosing $A(x)=0$ and $B(x)=I$ in (\ref{ljung_rec_elso}). Although this is the subject of the classic paper Robbins-Monroe-scheme, see \cite{robbins-monroe},
to have a unified treatment we shall discuss this problem using the DFL-scheme.
A possible motivation of this is the problem of identifying the noise characteristics of a L\'evy process using i.i.d. samples $y_1,y_2,\ldots$ generated by the increments of the process.

We suppose that the characteristic function of $y_i$ is known up to an unknown parameter $\eta^*.$ Let the c.f. of $y_i$ denoted by $\varphi(u,\eta^*).$ Fix a set of real $u_i$-s $1 \leq i \leq M.$
In this case following the idea of the off-line ECF method
our aim is to solve the non-linear equation $F(x)=0$ in (\ref{eq:G=0}) with $x=\eta$ and
$$
F(\eta)=\E \left[ -\varphi^*_{\eta}(\eta)K^{-1}h_N(\eta)\right]=0,
$$
where
$$\varphi_{\eta}(\eta)=\left(\varphi_{\eta}(u_1,\eta),\ldots,\varphi_{\eta}(u_M,\eta)\right)^T$$
and
\begin{equation}
\begin{split}
h_N(\eta)= \left( e^{iu_1 y_N}-\varphi(u_1,\eta),\ldots,e^{iu_M y_N}-\varphi(u_M,\eta)\right)^T.  
\end{split}
\end{equation}
A stochastic Newton method corresponding to this equation would read as
$$
\hat{\eta}_{N-}=\hat{\eta}_{N-1}-\frac1N (R_E^*)^{-1} \left( -\hat{\varphi}^*_{\eta,N}K^{-1}\hat{h}_N\right),
$$
with $R_E^*=\varphi^*_\eta(\eta^*)K^{-1}\varphi_\eta(\eta^*).$ Since $R^*$ is unknown we estimate it using the most current estimate of $\eta^*.$
Hence we extend the parameter vector $\eta$ to $(\eta,R)$ and re-define the equation $F(\eta)=0$ as
$$
F(\eta,R)=
\E\left[
\begin{array}{c}
-\varphi^*_{\eta}(\eta)K^{-1}h_N(\eta)\\
\varphi^*_{\eta}(\eta)K^{-1}\varphi_{\eta}(\eta)-R\\
\end{array}
\right]=
\left(
\begin{array}{c}
0\\
0\\
\end{array}
\right).
$$
Let the variables $\hat{\varphi}_{\eta,N}$ and $\hat{h}_N$ be obtained by using the most current estimate of $\eta^*,$ that is
\begin{equation}
\begin{aligned}
\hat{\varphi}_{\eta,N}&=\varphi_{\eta,N}(\hat{\eta}_{N-1})=\left(\varphi_{\eta}(u_1,\hat{\eta}_{N-1}),\ldots,\varphi_{\eta}(u_M,\hat{\eta}_{N-1})\right)^T \\
\hat{h}_N&=h_N(\hat{\eta}_{N-1})=\left( e^{iu_1 y_n}-\varphi(u_1,\hat{\eta}_{N-1}),\ldots,e^{iu_M y_n}-\varphi(u_M,\hat{\eta}_{N-1}) \right)^T.
\end{aligned}
\end{equation}
Let $\hat{\eta}_0$ and $\hat{R}_0$ be initial guesses and let $\hat{\eta}_N$ and $\hat{R}_N$ be computed using a partially stochastic Newton method as follows:
{\bf Algorithm 1}[Recursive i.i.d. ECF method]

\begin{equation}\label{eq:recursion_eta}
\begin{aligned}
\hat{\eta}_{N-}&=\hat{\eta}_{N-1}-\frac1N \hat{R}^{-1}_{N-1} \left( -\hat{\varphi}^*_{\eta,N}K^{-1}\hat{h}_N\right) \\
\hat{R}_{N-}&=\hat{R}_{N-1}+\frac1N \left(\hat{\varphi}^*_{\eta,N} K^{-1} \hat{\varphi}_{\eta,N}- \hat{R}_{N-1}\right).
\end{aligned}
\end{equation}
We note in passing that instead of applying a recursion for computing $\hat{R}_{N-}$ we could define it by simple substitution:
$$
\hat{R}_{N-}=\varphi^*_\eta(\hat{\eta}_{N})K^{-1}\varphi_\eta(\hat{\eta}_{N}).
$$
The reason behind our choice is that it fits into the general framework of DFL-scheme.
Denote the expected value of $h_N(\eta)$ by $g(\eta)$:
$$
g(\eta)=\E\left[h_N(\eta)\right]=\left( \varphi(u_1,\eta^*)-\varphi(u_1,\eta),\varphi(u_M,\eta^*)-\varphi(u_M,\eta) \right)^T.
$$
The corresponding associated ODE with the extended variable is given by
\begin{equation}\label{eq:iid_ode}
\begin{aligned}
\dot{\eta}_t&=-R^{-1}_t \left( -\varphi^*_\eta(\eta_t)K^{-1}g(\eta_t)\right) \\
\dot{R_t}&=\varphi^*_\eta(\eta_t)K^{-1}\varphi_\eta(\eta_t)-R_t
\end{aligned}
\end{equation}
for $t>0.$
We have seen that in case of an off-line identification method the optimal choice of $K$ is $K=C.$
Recall the notation $C$ with entries
$$
C_{k,l} = \varphi(u_k - u_l,\eta^*) - \varphi(u_k ,\eta^*)
\varphi(- u_l,\eta^*).
$$
Obviously for i.i.d. samples most of the conditions of Theorem \ref{thm:rateconv_rec} are satisfied: in (\ref{ljung_rec_elso}) we have $A=0$ and $B=I,$ furthermore we have that $(y_n)$ is strictly stationary and $L$-mixing.
Therefore Theorem \ref{thm:rateconv_rec} and Corollary 1 imply the following result:
\begin{theorem}
Suppose that we are given an i.i.d. data generated by a random variable $X$ such that $X^2$ is in $M^*,$ suppose further that the differential equation (\ref{eq:iid_ode}) satisfies Condition \ref{cond:od_rate}. Then the Algorithm 1 equipped with resetting is convergent and we have
$$\hat{\eta}_N-\eta^*=O_M(N^{-1/2}).$$
Furthermore, if $K=C$ then the asymptotic covariance matrix $\Sigma^{{\rm (rec)}}_{\eta \eta}$ of $\hat{\eta}_N$ is
given by
$$
\Sigma^{{\rm (rec)}}_{\eta \eta} = (\varphi^*_{\eta}(\eta^*) C^{-1} \varphi_{\eta}(\eta^*))^{-1}.
$$
Hence, the estimate $\hat{\eta}_N$ is essentially asymptotically efficient.
\end{theorem}
Note that the local stability of the ODE with $\alpha > 1/2$ follows. For, the Jacobian of the ordinary differential equation at $\eta=\eta^*$ and $R=R^*$ is
$$
\left(
  \begin{array}{cc}
    -I & 0 \\
    J_{2,1} & -I \\
  \end{array}
\right),
$$
thus each eigenvalue of the Jacobian is -1, hence the top Lyapunov exponent can be chosen to be equal to $-1+c$ with any $c>0,$ which implies that the ODE is locally stable with $\alpha > 1/2.$ The structure of the above Jacobian also shows that the proposed method is a partially stochastic Newton method w.r.t. $\eta.$

\section{Recursive ECF for linear L\'evy systems with known noise characteristics}\label{sec:rec_arma}

This section is devoted to the presentation of a recursive ECF method for linear L\'evy systems with known noise characteristics. The definition of linear L\'evy systems is given by:
\begin{equation}
\label{eq:disc_levy2} \Delta y=A(\theta^*, q^{-1}) \Delta L,
\end{equation}
defined for the time range $- \infty < n < + \infty,$ where
$\Delta L_n$ is the increment of a L\'evy process $(L_t)$ with $
-\infty < t < + \infty$, and $L_0=0$,
%
%
%
over an interval $[(n-1)h,nh),$ with $h>0$ being a fixed sampling
interval, and $ -\infty < n < + \infty$.
Let us assume that a state space representation in
innovation form for this model is given by
\begin{align}
\Delta X_{n+1}&=H(\theta^*) \Delta X_n + K(\theta^*) \Delta L_n \\
\Delta Y_n&=T(\theta^*) \Delta X_n + \Delta L_n.
\end{align}
Then $H(\theta)-K(\theta)T(\theta)$ is the state transition matrix of the inverse process.
We will need the following stability conditions:
\begin{cond} \label{cond:stability} It is assumed that the system matrix $H(\theta^*)$ is stable and
$H(\theta) - K(\theta)T(\theta)$ are jointly stable for $\theta \in
D_{\theta}.$
\end{cond}
The next condition guarantees the smooth dependence of the system matrices on $\theta:$
\begin{cond}\label{cond:A_smooth}
Assume that
$H(\theta),K(\theta)$ and $T(\theta)$ are three-times continuously differentiable
w.r.t. $\theta$ for $\theta \in D_{\theta}$.
\end{cond}
The novel problem of identifying the system
parameters $\theta^*,$ using the ECF
method, under the assumption that the noise characteristics
$\eta^*$ is known was presented in \cite{automatica_own}.
We may wish to solve the same problem, but now with a recursive method.

Suppose that we are given the noise characteristics $\eta^*.$
Fix a set of real $u_i$-s $1 \leq i \leq M.$
Following the third step of the off-line estimation method we seek the solution of the non-linear
equation $F(x)=0$ with $x=\theta$ and
$$
F(\theta)=\E \left[ G^*K^{-1}h^{{\rm (s)}}_N(\theta;\eta^*)\right],
$$
where
\begin{equation*}
h^{{\rm (s)}}_N(\theta;\eta^*)=\left(\left( e^{iu_1 \varepsilon^{{\rm (s)}}_{N}(\theta)}-\varphi(u_1,\eta^*)\right)\varepsilon^{{\rm (s)T}}_{\theta N}(\theta),\ldots,\left(e^{iu_M \varepsilon^{{\rm (s)}}_{N}(\theta)}-\varphi(u_M,\eta^*)\right)\varepsilon^{{\rm (s)T}}_{\theta N}(\theta)\right)^T.
\end{equation*}
Clearly $\theta=\theta^*$ is the solution of this equation.
Similarly to the i.i.d. case we would like to apply a stochastic Newton method, which requires the introduction of the Jacobian $R^*,$ the true value of which is known to be $R^*=G^*(\theta^*,\eta^*)K^{-1}G(\theta^*,\eta^*)$. In order to compute $R^*$ we need the value of $G(\theta^*,\eta^*)=\psi \otimes R_P^*,$ where $\psi=\left( iu_1 \varphi(u_1,\eta^*),\ldots,iu_M \varphi(u_M,\eta^*)\right)^T,$ for $G(\theta^*,\eta^*)$ see the proof of Theorem 8 in \cite{automatica_own}. Since these are not computable without the knowledge of $\theta^*,$ and $G(\theta^*,\eta^*)$ can be computed only empirically, we approximate $G(\theta^*,\eta^*)$ using the most current estimate of $\theta^*.$
To this end we extend the parameter vector to $(\theta,G,R)$ and re-define equation $F(\theta)=0$ by
\begin{equation}\label{eq:F3}
F(\theta,G,R)=
\E\left[
\begin{array}{c}
G^*K^{-1}h^{{\rm (s)}}_N(\theta;\eta)\\
G-h^{{\rm (s)}}_{\theta N}(\theta;\eta)\\
G^*K^{-1}G-R \\
\end{array}
\right],
\end{equation}
where $h^{{\rm (s)}}_{\theta N}(\theta)$ shows up in the derivative of $h^{{\rm (s)}}_{N}(\theta)$ w.r.t. $\theta,$ and defined by
\begin{equation}
h^{{\rm (s)}}_{\theta N}(\theta):=\left( \left(i u_1 e^{iu_1\varepsilon^{{\rm (s)}}_N(\theta)}\varepsilon^{{\rm (s)}}_{\theta,N}(\theta)\varepsilon^{{\rm (s)*}}_{\theta,N}(\theta)\right)^T,\ldots,\left(i u_M e^{iu_M \varepsilon^{{\rm (s)}}_N(\theta)}\varepsilon^{{\rm (s)}}_{\theta,N}(\theta)\varepsilon^{{\rm (s)*}}_{\theta,N}(\theta)\right)^T\right)^T,
\end{equation}
which is obtained by dropping the term containing $-u^2_j e^{iu_j\varepsilon^{{\rm (s)}}_N(\theta)}\varepsilon^{{\rm (s)}}_{\theta \theta,N}(\theta)$-s,
which has zero expectation at $\theta=\theta^*,$ from the derivative of $h^{{\rm (s)}}_{N}(\theta)$.

Suppose we are given a set of initial values of the parameter: $\hat{\theta}_0 \in D_{0 \theta}$ is the initial value of $\theta$ and $\hat{g}_{\theta,0}$ is the initial value of $G$. The set of initial values of the parameter is given by $\left\{ \hat{\varepsilon}_0, \hat{\varepsilon}_{\theta,0} \right\}.$ Here $\hat{\theta}_0$ might have been previously obtained by an off-line identification method. Likewise $\hat{\varepsilon}_0, \hat{\varepsilon}_{\theta,0},\hat{g}_{\theta,0}$ might be obtained previously or we can set them to be equal to 0.
The recursive algorithm at step $N$ updates the estimates as follows: given the previous estimates first compute the auxiliary variables
using the most current estimate of $\theta^*$ according to
\begin{equation}
\begin{aligned}
\hat{\varepsilon}_N&=A^{-1}\left(\hat{\theta}_{N-1}\right)\Delta y_N \\
\hat{\varepsilon}_{\theta,N}&=A_{\theta}^{-1}\left(\hat{\theta}_{N-1}\right)\Delta y_N \\
\hat{h}_N(\eta^*)&=\left(\left( e^{iu_1 \hat{\varepsilon}_N}-\varphi(u_1,\eta^*)\right)\hat{\varepsilon}^T_{\theta,N},\ldots,\left( e^{iu_M \hat{\varepsilon}_N}-\varphi(u_M,\eta^*)\right)\hat{\varepsilon}^T_{\theta,N}\right)^T\\
\hat{h}_{\theta,N}&=\left( \left(i u_1 e^{iu_1 \hat{\varepsilon}_N}\hat{\varepsilon}_{\theta,N}\hat{\varepsilon}^*_{\theta,N}\right)^T,\ldots,\left(i u_M e^{iu_M \hat{\varepsilon}_N}\hat{\varepsilon}_{\theta,N}\hat{\varepsilon}^*_{\theta,N}\right)^T\right)^T.
\end{aligned}
\end{equation}
Following the special form of the off-line estimation method presented in \cite{automatica_own} we define a stochastic Newton method via the following algorithm:

{\bf Algorithm 2}[Re-estimating recursive ECF method]

\begin{equation}\label{eq:recursion_theta}
\begin{aligned}
\hat{\theta}_{N-}&=\hat{\theta}_{N-1}-\frac1N \hat{R}_{S,N-1}^{-1} \left( \hat{g}^*_{\theta,N-1}K^{-1}\hat{h}_N(\eta^*)\right) \\
\hat{g}_{\theta,N-}&=\hat{g}_{\theta,N-1}+\frac1N \left( \hat{h}_{\theta,N}-\hat{g}_{\theta,N-1}\right),
\end{aligned}
\end{equation}
where
$$\hat{R}_{S,N-1}=\hat{g}^*_{\theta,N-1}K^{-1}\hat{g}_{\theta,N-1}.$$
The Jacobian of the
Note that the third component $G^*K^{-1}G-R$ in (\ref{eq:F3}) is non-random, hence $\hat{R}_{S,N-1}$ is computed by simple substitution.
These tentative values need to be modified with a suitable resetting mechanism as described in connection with the general DFL-scheme.
In order to define the associated ODE first take the expectations of the frozen parameter correction terms $h^{{\rm (s)}}_{\theta N}(\theta;\eta)$ and $h^{{\rm (s)}}_{N}(\theta;\eta),$ showing up on the r.h.s. of \ref{eq:F3}:
\begin{align*}
g(\theta;\eta^*)&=\E\left[\left(\left(e^{iu_1\varepsilon^{{\rm (s)}}_n(\theta)}-\varphi(u_1,\eta^*)\right)\varepsilon^{{\rm (s)T}}_{\theta,n}(\theta),\ldots,\left(e^{iu_M\varepsilon^{{\rm (s)}}_n(\theta)}-\varphi(u_M,\eta^*)\right)\varepsilon^{{\rm (s)T}}_{\theta,n}(\theta)\right)^T\right],\\
g_{\theta}(\theta)&=\E\left[\left( \left(i u_1 e^{iu_1\varepsilon^{{\rm (s)}}_n(\theta)}\varepsilon^{{\rm (s)}}_{\theta,n}(\theta)\varepsilon^{{\rm (s)*}}_{\theta,n}(\theta)\right)^T,\ldots,\left(i u_M e^{iu_M \varepsilon^{{\rm (s)}}_n(\theta)}\varepsilon^{{\rm (s)}}_{\theta,n}(\theta)\varepsilon^{{\rm (s)*}}_{\theta,n}(\theta)\right)^T\right)^T \right].
\end{align*}
Then the corresponding associated ODE reads as
\begin{equation}\label{eq:2steps_ode}
\begin{aligned}
\dot{\theta}_{t}&=-R_{S,t}^{-1}\left( g^*_{\theta,t}K^{-1}g(\theta_{t};\eta^*)\right) \\
\dot{g}_{\theta,t}&=g_{\theta}(\theta_{t})-g_{\theta,t},
\end{aligned}
\end{equation}
where
$$R_{S,t}=g^*_{\theta,t}K^{-1}g_{\theta,t}.$$
It is easy to check that the Jacobian matrix of (\ref{eq:2steps_ode}) at $(\theta^*,G(\theta^*,\eta^*))$ is a lower triangular matrix with $-I$ blocks in the diagonal, thus all eigenvalues are equal to $-1.$ It follows that the solution of the ODE is locally stable with $\alpha=1/2,$ see Condition \ref{cond:od_rate} for the definition of $\alpha$. Moreover, the structure of the Jacobian also shoves that the on-line method is a partially stochastic Newton method w.r.t. $\theta.$

Recall the notation of \cite{automatica_own}:
$$ R_{P}^*=\E \left[ \varepsilon^{(\rm s)}_{\theta n} (\theta^*)\varepsilon^{(\rm s)T}_{\theta n}(\theta^*) \right].$$
Theorem \ref{thm:rateconv_rec} and Corollary 1 together with Theorem 9 in \cite{automatica_own}, giving the asymptotic covariance matrix of the re-estimated system parameter, imply the following result:
\begin{theorem}\label{thm:error_3rd_step}
Suppose that for the L\'evy system Condition \ref{cond:stability} and Condition \ref{cond:A_smooth} are satisfied. Suppose further that for the driving L\'evy process we have that $((\Delta L_n)^2)$ is in $M^*$ and that the differential equation (\ref{eq:2steps_ode}) satisfies Condition \ref{cond:od_rate}.
Then for the estimate $\hat{\theta}_N$ obtained by the above recursive method in Algorithm 2 modified by a suitable resetting mechanism we have
$$ \hat{\theta}_N-\theta^*=O_M(N^{-1/2}). 
$$
Moreover, if $K=C \otimes R_{P}^*,$ then the asymptotic covariance matrix $\Sigma^{{\rm (rec)}}_{\theta \theta}$ of $\hat{\theta}_N-\theta^*$ is given by
$$
\Sigma^{{\rm (rec)}}_{\theta \theta}=\left(\psi^* C^{-1} \psi\right)^{-1} \left(R_P^*\right)^{-1},
$$
with $\psi=\left( iu_1 \varphi(u_1,\eta^*),\ldots,iu_M \varphi(u_M,\eta^*)\right)^T.$ Hence, the proposed on-line method is essentially asymptotically efficient.
\end{theorem}

\section{Recursive ECF for linear L\'evy systems}\label{sec:rec_3step}

Now we are ready to present and analyze a recursive identification method that estimates both the system and the noise characteristics by converting the three-stage method presented in \cite{automatica_own} to a recursive method. Suppose that the dynamics of $(y_n)$ follows (\ref{eq:disc_levy2}).
Fix a set of real $u_i$-s $1 \leq i \leq M.$
In this case we define the non-linear equation $F(x)=0$ in (\ref{eq:G=0}) by merging the asymptotic equations corresponding to the PE method, the ECF method for the noise characteristics $\eta^*$ and the ECF method for the re-estimation of $\theta^*.$ Accordingly, $x$ is defined as $x=(\theta_P,R_P,\eta,R_E,\theta_S,G,R_S).$ Observe that $\theta$ is duplicated, in the sense that $\hat{\theta}_P$ and $\hat{\theta}_S$ both are expected to converge to $\theta^*.$ This separation of the recursive PE estimate $\hat{\theta}_P$ and the ECF estimate $\hat{\theta}_S$ guarantees that the Jacobian matrix of the corresponding associated ODE will be lower triangular. Now we are ready to define $F$ by
$$
F(\theta_P,R_P,\eta,R_E,\theta_S,G,R_S)=
\E\left[
\begin{array}{c}
\varepsilon^{{\rm (s)}}_{\theta_P N}(\theta_P)\varepsilon^{{\rm (s)}}_{N}(\theta_P)\\
\varepsilon^{{\rm (s)}}_{\theta_P N}(\theta_P)\varepsilon^{{\rm (s)T}}_{\theta_P N}(\theta_P)-R_P\\ 
-\varphi^*_{\eta}(\eta)K_{E}^{-1}h^{{\rm (s)}}_{E, N}(\theta_P,\eta)\\
\varphi^*_{\eta}(\eta)K_{E}^{-1}\varphi_{\eta}(\eta)-R_E\\
G^*K_S^{-1}h^{{\rm (s)}}_{S, N}(\theta_S,\eta)\\
G-h^{{\rm (s)}}_{S, \theta N}(\theta_S,\eta)\\
G^*K_S^{-1}G-R_S \\
\end{array}
\right],
$$
where the auxiliary variables are defined in analogy with the ones in the previous two sections as
\begin{align*}
\varphi_{\eta}(\eta)&=\left(\varphi_{\eta}(u_1,\eta),\ldots,\varphi_{\eta}(u_M,\eta)\right)^T
\\
h^{{\rm (s)}}_{E,N}(\theta,\eta)&=\left( e^{iu_1 \varepsilon^{{\rm (s)}}_N(\theta)}-\varphi(u_1,\eta),\ldots,e^{iu_M \varepsilon^{{\rm (s)}}_N(\theta)}-\varphi(u_M,\eta)\right)^T
\\
h^{{\rm (s)}}_{S,N}(\theta,\eta)&=\left(\left( e^{iu_1 \varepsilon^{{\rm (s)}}_{N}(\theta)}-\varphi(u_1,\eta)\right)\varepsilon^{{\rm (s)T}}_{\theta N}(\theta),\ldots,\left(e^{iu_M \varepsilon^{{\rm (s)}}_{N}(\theta)}-\varphi(u_M,\eta)\right)\varepsilon^{{\rm (s)T}}_{\theta N}(\theta)\right)^T
\\
h^{{\rm (s)}}_{S,\theta N}(\theta)&=\left( \left(i u_1 e^{iu_1\varepsilon^{{\rm (s)}}_n(\theta)}\varepsilon^{{\rm (s)}}_{\theta,n}(\theta)\varepsilon^{{\rm (s)*}}_{\theta,n}(\theta)\right)^T,\ldots,\left(i u_M e^{iu_M \varepsilon^{{\rm (s)}}_n(\theta)}\varepsilon^{{\rm (s)}}_{\theta,n}(\theta)\varepsilon^{{\rm (s)*}}_{\theta,n}(\theta)\right)^T\right)^T.
\end{align*}
Note that to different ECF scores $h_E$ and $h_S$ are being used, one estimates $\eta^*$ and another re-estimates $\theta^*.$
Let us suppose that we are given the initial values of the parameters: $\hat{\theta}_{P,0},\hat{R}_{P,0}$ are the initial values of the recursive PE method, see \cite{soderstrom_rec}, $\hat{\eta}_{0},\hat{R}_{E,0}$ are the initial values of the recursive ECF method for the noise characteristics, see Section \ref{sec:rec_noise} and $\hat{\theta}_{S,0},\hat{g}_{\theta,0}$ are the initial values of the recursive ECF re-estimation method, see Section \ref{sec:rec_arma}. We assume that each of these initial values are the element of the corresponding truncation domain and $\hat{\theta}_{P,0}=\hat{\theta}_{S,0}$ is a reasonable choice. We are also given a set of initial values $\hat{\varepsilon}_{P,0}, \hat{\varepsilon}_{P,\theta,0},\hat{\varepsilon}_{S,0}, \hat{\varepsilon}_{S,\theta,0}.$ Clearly, these values might have been obtained by carrying out an off-line identification method, otherwise we may set all of them to be equal to zero.

The recursive algorithm at step $N$ updates the estimates as follows: given the previous estimates of the parameters first compute the estimated driving noise and its derivative using the most current values of the parameters as
\begin{equation}
\begin{aligned}
\hat{\varepsilon}_{P,N}&=A^{-1}\left(\hat{\theta}_{P,N-1}\right)\Delta y_N \\
\hat{\varepsilon}_{P,\theta,N}&=A_{\theta}^{-1}\left(\hat{\theta}_{P,\theta,N-1}\right)\Delta y_N \\
\hat{\varepsilon}_{S,N}&=A^{-1}\left(\hat{\theta}_{S,N-1}\right)\Delta y_N \\
\hat{\varepsilon}_{S,\theta,N}&=A_{\theta}^{-1}\left(\hat{\theta}_{S,\theta,N-1}\right)\Delta y_N.
\end{aligned}
\end{equation}
While the first two equations correspond to the recursive PE method, the last two equations correspond to the re-estimating ECF method for linear systems.
In analogy with the previous two sections we also define the auxiliary variables using the most current values of the parameters according to
\begin{equation}
\begin{aligned}
\hat{\varphi}_{\eta,N}&=\left(\varphi_{\eta}(u_1,\hat{\eta}_{N-1}),\ldots,\varphi_{\eta}(u_M,\hat{\eta}_{N-1})\right)^T \\
\hat{h}_{E,N}&=\left( e^{iu_1 \hat{\varepsilon}_{P,N}}-\varphi(u_1,\hat{\eta}_{N-1}),\ldots,e^{iu_M \hat{\varepsilon}_{P,N}}-\varphi(u_M,\hat{\eta}_{N-1}) \right)^T \\
\hat{h}_{S,N}&=\left(\left( e^{iu_1 \hat{\varepsilon}_{S,N}}-\varphi(u_1,\hat{\eta}_{N-1})\right)\hat{\varepsilon}^T_{S,\theta,N},\ldots,\left( e^{iu_M \hat{\varepsilon}_{S,N}}-\varphi(u_M,\hat{\eta}_{N-1})\right)\hat{\varepsilon}^T_{S,\theta,N}\right)^T\\
\hat{h}_{S,\theta,N}&=\left( \left(i u_1 e^{iu_1 \hat{\varepsilon}_{S,N}}\hat{\varepsilon}_{S,\theta,N}\hat{\varepsilon}^*_{S,\theta,N}\right)^T,\ldots,\left(i u_M e^{iu_M \hat{\varepsilon}_{S,N}}\hat{\varepsilon}_{S,\theta,N}\hat{\varepsilon}^*_{S,\theta,N}\right)^T\right)^T.
\end{aligned}
\end{equation}
The recursive version of the three-stage method is then given as follows:

{\bf Algorithm 3}[Three-stage recursive ECF method]

First apply the recursive PE method defined as
\begin{equation}\label{eq:recursion_thetaeta_1}
\begin{aligned}
\hat{\theta}_{P,N-}&=\hat{\theta}_{N-1}-\frac1N \hat{R}_{P,N-1}^{-1} \hat{\varepsilon}_{P,\theta N} \hat{\varepsilon}^T_{P,N} \\
\hat{R}_{P,N-}&=\hat{R}_{E,N-1}+\frac1N \left(\hat{\varepsilon}_{P,\theta N}\hat{\varepsilon}^T_{P,\theta N}- \hat{R}_{P,N-1}\right),
\end{aligned}
\end{equation}
then apply the recursive ECF method for the noise characteristics defined via
\begin{equation}\label{eq:recursion_thetaeta_2}
\begin{aligned}
\hat{\eta}_{N-}&=\hat{\eta}_{N-1}-\frac1N \hat{R}_{E,N-1}^{-1} \left( -\hat{\varphi}^*_{\eta,N}K_E^{-1}\hat{h}_{E,N}\right) \\
\hat{R}_{E,N-}&=\hat{R}_{E,N-1}+\frac1N \left(\hat{\varphi}^*_{\eta,N} K_E^{-1} \hat{\varphi}_{\eta,N}- \hat{R}_{P,N-1}\right),
\end{aligned}
\end{equation}
finally re-estimate $\theta^*$ using the recursive ECF method defined by
\begin{equation}\label{eq:recursion_thetaeta_3}
\begin{aligned}
\hat{\theta}_{S,N-}&=\hat{\theta}_{N-1}-\frac1N \hat{R}_{S,N-1}^{-1} \left( \hat{g}^*_{\theta,N-1}K_S^{-1}\hat{h}_{S,N}\right) \\
\hat{g}_{\theta,N-}&=\hat{g}_{\theta,N-1}+\frac1N \left( \hat{h}_{S,\theta,N}-\hat{g}_{\theta,N-1}\right),
\end{aligned}
\end{equation}
where $\hat{R}_{S,N-1}=\hat{g}^*_{\theta,N-1}K^{-1}\hat{g}_{\theta,N-1}$.

These tentative values need to be modified using a suitable resetting mechanism as described in connection with the DFL-scheme.
Write the expectations of the frozen parameters as
\begin{align*}
R_P(\theta_P)&=\E\left[ \varepsilon^{{\rm (s)}}_{\theta,n}(\theta_P)\varepsilon^{{\rm (s)T}}_{\theta,n}(\theta_P)\right] \\
h_E(\theta_P,\eta)&=\E\left[ \left( e^{iu_1\varepsilon^{{\rm (s)}}_{n}(\theta_P)}-\varphi(u_1,\eta),\ldots, e^{iu_M\varepsilon^{{\rm (s)}}_{n}(\theta_P)}-\varphi(u_M,\eta)\right)^T\right] \\
h_S(\theta_S,\eta)&=\E\left[\left(\left(e^{iu_1\varepsilon^{{\rm (s)}}_n(\theta_S)}-\varphi(u_1,\eta)\right)\varepsilon^{{\rm (s)T}}_{\theta,n}(\theta_S),\ldots,\left(e^{iu_M\varepsilon^{{\rm (s)}}_n(\theta_S)}-\varphi(u_M,\eta)\right)\varepsilon^{{\rm (s)T}}_{\theta,n}(\theta_S)\right)^T\right]\\
g_{\theta}(\theta_S)&=\E\left[\left( \left(i u_1 e^{iu_1\varepsilon^{{\rm (s)}}_n(\theta_S)}\varepsilon^{{\rm (s)}}_{\theta,n}(\theta_S)\varepsilon^{{\rm (s)*}}_{\theta,n}(\theta_S)\right)^T,\ldots,\left(i u_M e^{iu_M \varepsilon^{{\rm (s)}}_n(\theta_S)}\varepsilon^{{\rm (s)}}_{\theta,n}(\theta_S)\varepsilon^{{\rm (s)*}}_{\theta,n}(\theta_S)\right)^T\right)^T \right].
\end{align*}
Recall the notation
$$
W_{P,\theta_P}(\theta_{P})=\E \left[ \varepsilon^{{\rm (s)}}_{\theta n}(\theta_P)\varepsilon^{{\rm (s)}}_n(\theta_P)\right].
$$
In terms of the above expectations the ODE corresponding to the recursive PE method reads as
\begin{equation}\label{eq:3steps_ode_1}
\begin{aligned}
\dot{\theta}_{P,t}&=-R_{P,t}^{-1}W_{P,\theta_P}(\theta_{P,t}) \\
\dot{R}_{P,t}     &=R_P(\theta_{P,t})-R_{P,t},
\end{aligned}
\end{equation}
while that of the recursive ECF method for noise characteristic is given by
\begin{equation}\label{eq:3steps_ode_2}
\begin{aligned}
\dot{\eta}_t      &=-R^{-1}_{E,t} \left( -\varphi^*_\eta(\eta_t)K^{-1}h_E(\theta_{P,t},\eta_t)\right) \\
\dot{R}_{E,t}     &=\varphi^*_{\eta}(\eta_t)K_E^{-1}\varphi_\eta(\eta_t)-R_{E,t}
\end{aligned}
\end{equation}
and finally the ODE of the ECF method system parameters can be written as
\begin{equation}\label{eq:3steps_ode_3}
\begin{aligned}
\dot{\theta}_{S,t}&=-R_{S,t}^{-1}\left( g^*_{\theta,t}K_S^{-1}h_S(\theta_{S,t},\eta_t)\right) \\
\dot{g}_{\theta,t}&=g_{\theta}(\theta_{S,t})-g_{\theta,t},
\end{aligned}
\end{equation}
where $R_{S,t}=g^*_{\theta,t}K^{-1}g_{\theta,t}.$
By merging the above three ODE-s we get the associated ODE of the recursive three-stage identification method:
\begin{equation}\label{eq:3steps_ode}
\begin{aligned}
\dot{\theta}_{P,t}&=-R_{P,t}^{-1}W_{P,\theta_P}(\theta_{P,t}) \\
\dot{R}_{P,t}     &=R(\theta_{P,t})-R_{P,t} \\
\dot{\eta}_t      &=-R^{-1}_{E,t} \left( -\varphi^*_\eta(\eta_t)K^{-1}h_E(\theta_{P,t},\eta_t)\right) \\
\dot{R}_{E,t}     &=\varphi^*_{\eta}(\eta_t)K_E^{-1}\varphi_\eta(\eta_t)-R_{E,t} \\
\dot{\theta}_{S,t}&=-R_{S,t}^{-1}\left( g^*_{\theta,t}K_S^{-1}h_S(\theta_{S,t},\eta_t)\right) \\
\dot{g}_{\theta,t}&=g_{\theta}(\theta_{S,t})-g_{\theta,t}.
\end{aligned}
\end{equation}
The Jacobian of the r.h.s. at
$$
(\theta_P,R_P,\eta,R_E,\theta_S,G)=(\theta^*,R_P^*,\eta^*,R_E^*,\theta^*,G(\theta^*,\eta^*))
$$
is given by
$$
\left(
  \begin{array}{cccccc}
    -I & 0       & 0 & 0 & 0 & 0 \\
    J_{2,1} & -I & 0 & 0 & 0 & 0 \\
    0 & 0        & -I & 0 & 0 & 0 \\
    0 & 0        & J_{4,3} & -I & 0 & 0 \\
    0 & 0        & 0       & 0 & -I & 0 \\
    0 & 0        & 0       & 0 & J_{6,5}  & -I\\
  \end{array}
\right).
$$
Hence the solution of the ODE \ref{eq:3steps_ode} is locally stable with $\alpha=1/2,$ see Condition \ref{cond:od_rate} for the definition of $\alpha$.
The structure of the above Jacobian matrix and Theorem \ref{thm:rateconv_rec} together with Corollary 1 imply the next result.
\begin{theorem}\label{thm:error}
Let $\hat{\theta}_{S,N}$ and $\hat{\eta}_N$ be the $N^{th}$-step estimate of the parameters obtained by the recursive estimation in Algorithm 3
 using a suitable resetting mechanism. Suppose that for the L\'evy system Condition \ref{cond:stability} and Condition \ref{cond:A_smooth} are satisfied. Suppose further that for the driving L\'evy process we have that $((\Delta L_n)^2)$ is in $M^*$ and that the differential equation (\ref{eq:3steps_ode}) satisfies Condition \ref{cond:od_rate}. 
Then
we have
$$ \hat{\eta}_N-\eta^*=O_M(N^{-1/2}) \text{~~and~~} \hat{\theta}_{S,N}-\theta^*=O_M(N^{-1/2}).$$
Furthermore, if $K_E=C$ then the asymptotic covariance matrix $\Sigma^{{\rm (rec)}}_{\eta \eta}$ of $\hat{\eta}_N$ is given by
$$
\Sigma^{{\rm (rec)}}_{\eta \eta} = (\varphi^*_{\eta}(\eta^*) C^{-1} \varphi_{\eta}(\eta^*))^{-1},
$$
and if $K_S=C \otimes R_{P}^*,$ then the asymptotic covariance matrix $\Sigma^{{\rm (rec)}}_{\theta \theta}$ of $\hat{\theta}_{S,N}-\theta^*$ is given by
$$
\Sigma^{{\rm (rec)}}_{\theta \theta}=\left(\psi^* C^{-1} \psi\right)^{-1} \left(R_P^*\right)^{-1},
$$
with $\psi=\left( iu_1 \varphi(u_1,\eta^*),\ldots,iu_M \varphi(u_M,\eta^*)\right)^T.$ Hence, the proposed on-line method gives essentially asymptotically efficient estimates of $\theta^*$ and $\eta^*.$
\end{theorem}

\emph{Remark:}

Similarly to the off-line identification the optimal choice of $K_E$ and $K_S$ depend on the true values $\theta^*, \eta^*.$ Thus these values need to be approximated using the most recent estimates of the parameters. In any case it can be easily shown that the results of this paper remain valid even if these approximated weighting matrices are used.


\bibliography{refereces_new}{}
\bibliographystyle{abbrv}
\end{document}